\documentclass[12pt]{amsart}
\usepackage{amssymb}
\usepackage{mathtools}
\usepackage[english]{babel}
\usepackage{epsfig}
\usepackage{amscd}
\usepackage{amsmath}
\usepackage{graphicx}
\usepackage{wasysym}
\usepackage{enumerate}
\usepackage{ifpdf}
\usepackage{amsfonts}
\usepackage{amsthm}
\usepackage[hyperfootnotes=false]{hyperref}
\hypersetup{hidelinks}
\usepackage{setspace}
\usepackage[usenames]{xcolor}
\usepackage{tikz}
\usepackage[normalem]{ulem}

\numberwithin{equation}{section}

\def\Re{{\sf Re}\,}
\def\Im{{\sf Im}\,}

\def\ech{.6}

\def\1#1{\overline{#1}}
\def\2#1{\widetilde{#1}}
\def\3#1{\widehat{#1}}
\def\4#1{\mathbb{#1}}
\def\5#1{\frak{#1}}
\def\6#1{{\mathcal{#1}}}

\newcommand{\C}{\mathbb C}

\newcommand{\D}{\mathbb D}
\newcommand{\oD}{\overline{\mathbb D}}

\newcommand{\N}{\mathbb N}

\def\Re{{\sf Re}\,}
\def\Im{{\sf Im}\,}

\newcommand{\mcite}[1]{\csname b@#1\endcsname}

\theoremstyle{theorem}

\def\Re{{\sf Re}\,}
\def\Im{{\sf Im}\,}

\newtheorem{theorem}{Theorem}[section]
\newtheorem{lemma}[theorem]{Lemma}
\newtheorem{proposition}[theorem]{Proposition}
\newtheorem{corollary}[theorem]{Corollary}

\theoremstyle{definition}
\newtheorem{definition}[theorem]{Definition}
\newtheorem{example}[theorem]{Example}

\theoremstyle{remark}
\newtheorem{remark}[theorem]{Remark}

\numberwithin{equation}{section}

\newcommand{\abs}[1]{\left|#1\right|}

\title[Asymptotic behavior of orbits]{Asymptotic behavior of orbits of holomorphic semigroups}

\author[F. Bracci]{Filippo Bracci$^1$}
\address{F. Bracci: Dipartimento di Matematica, Universit\`a di Roma ``Tor Vergata", Via della Ricerca
Scientifica 1, 00133, Roma, Italia.} \email{fbracci@mat.uniroma2.it}

\author[M. D. Contreras]{Manuel D. Contreras$^2$}

\author[S. D\'{\i}az-Madrigal]{Santiago D\'{\i}az-Madrigal$^{2}$}
\address{M. D. Contreras, S. D\'{\i}az-Madrigal: Camino de los Descubrimientos, s/n\\
Departamento de Matem\'{a}tica Aplicada~II and IMUS\\ Universidad de Sevilla\\ Sevilla,
41092\\ Spain.}\email{contreras@us.es} \email{madrigal@us.es}

\author[H. Gaussier]{Herv\'e Gaussier$^3$}
\address{H. Gaussier: Univ. Grenoble Alpes, CNRS, IF, F-38000 Grenoble, France} \email{herve.gaussier@univ-grenoble-alpes.fr}

\author[A. Zimmer]{Andrew Zimmer$^4$}
\address{A. Zimmer: Department of Mathematics, Louisiana State University, Baton Rouge, LA USA} \email{amzimmer@lsu.edu}

\date{\today}
\thanks{$^1\,$Partially supported by the MIUR Excellence Department Project awarded to the  
Department of Mathematics, University of Rome Tor Vergata, CUP E83C18000100006}
\thanks{$^{2}\,$Partially supported by the \textit{Ministerio
de Econom\'{\i}a y Competitividad} and the European Union (FEDER) MTM2015-63699-P and by \textit{La Consejer\'{\i}a de Educaci\'{o}n y Ciencia de la Junta de Andaluc\'{\i}a}.}
\thanks{$^3\,$Partially supported by ERC ALKAGE}
\thanks{$^4\,$Partially supported by the National Science Foundation under grant DMS-1760233.}

\long\def\REM#1{\relax}

\begin{document}

\begin{abstract}
Let $(\phi_t)$ be a holomorphic semigroup of the unit disc ({\sl i.e.}, the flow of a semicomplete holomorphic vector field) without fixed points in the unit disc and let $\Omega$ be the starlike at infinity domain image of the Koenigs function of $(\phi_t)$. In this paper we completely characterize the type of convergence of the orbits of $(\phi_t)$ to the Denjoy-Wolff point in terms of the shape of $\Omega$. In particular we prove that the convergence is non-tangential if and only if the domain $\Omega$ is ``quasi-symmetric with respect to vertical axes''. We also prove that such conditions are equivalent to the curve $[0,\infty)\ni t\mapsto \phi_t(z)$ being a  quasi-geodesic in the sense of Gromov. Also, we characterize  the tangential convergence in terms of the shape of $\Omega$. 
\end{abstract}

\maketitle


\section{Introduction and statements of the main results}

A holomorphic vector field $G$ on the unit disc $\D$ is (real) semicomplete if the Cauchy problem $\dot{x}(t)=G(x(t)), x(0)=z$ has a solution defined for all $t\geq 0$ and for all $z\in \D$. The flow of a semicomplete vector field, $(\phi_t)$, is a continuous semigroup of holomorphic self-maps of $\D$---or simply a  semigroup in $\D$. Namely, $(\phi_t)$ is a continuous homomorphism of the real semigroup $[0,+\infty)$ endowed with the Euclidean topology to the semigroup under composition of holomorphic self-maps of $\D$ endowed with the topology of uniform convergence on compacta.

It appears that semigroups in $\D$ were first considered in the 1930's by J. Wolff \cite{Wo}, although it was only with a paper of E. Berkson and H. Porta \cite{Berkson-Porta} in the 1970's that  the modern  study of semigroups in $\D$ initiated. Since their work, interest in semigroups in $\D$ has expanded due to various applications in physics and biology (see, {\sl e.g.},  \cite{KaMc1, KaMc2, Mit}) and their connections to composition operators and Loewner's theory (we refer the reader to the books \cite{Abate, Shb, ReiShobook05, EliShobook10} and \cite{BCD} for more details). 

In this paper we completely characterize the  asymptotic behavior of semigroups in $\D$ via the Euclidean geometry of the image of an associated Koenigs function. Aside being motived by the study of the dynamics of semigroups, our main results also give a complete answer to the following  question from geometric function theory: 

\smallskip

{\sl Let $f:\D\to \C$ be a Riemann map such that $\Omega:=f(\D)$ is starlike at infinity. Let $p\in \Omega$ and let $\{t_n\}$ be a sequence of positive real numbers converging to $+\infty$. Looking only at the shape of $\Omega$, how can one decide whether the sequence $\{f^{-1}(p+it_n)\}$ converges to a point $\tau\in \partial \D$ non-tangentially or tangentially?}

\smallskip

If $(\phi_t)$ is a semigroup in $\D$, which is not a group of hyperbolic rotations, then there exists a unique $\tau\in \oD$, the {\sl Denjoy-Wolff point} of $(\phi_t)$, such that $\lim_{t\to +\infty}\phi_t(z)=\tau$, and the convergence is uniform on compacta. In case $\tau\in \D$, the semigroup is called \emph{elliptic}. 

In case the semigroup $(\phi_t)$ is non-elliptic, the action is conjugate to linear translation on an unbounded simply connected domain. More precisely, there exists  an (essentially unique) univalent function $h$, called the {\sl Koenigs function} of $(\phi_t)$, such that $h(\D)$ is starlike at infinity, $h(\phi_t(z))=h(z)+it$ for all $t\geq 0$ and $z\in \D$ (see, {\sl e.g.}, \cite{Abate, BrAr, Cowen}).

The {\sl slope} of a non-elliptic semigroup $(\phi_t)$ at $z\in \D$ is the cluster set of ${\sf Arg}(1-\overline{\tau}\phi_t(z))$ as $t\to +\infty$. The slope is a compact connected subset of $[-\pi/2,\pi/2]$. 

Given $z\in \D$, we say that the orbit $[0,+\infty)\ni t\mapsto \phi_t(z)$ converges {\sl non-tangentially} to the Denjoy-Wolff point if the slope of $(\phi_t)$ at $z$ is contained in $(-\pi/2,\pi/2)$. In case the slope is $\{-\pi/2\}$ or $\{\pi/2\}$, the convergence is {\sl tangential}.

 For one-parameter groups of automorphisms there are two possible behaviors. Either  $h(\D)$ is a vertical strip (and the group is called {\sl hyperbolic}) or $h(\D)$ is a vertical half-plane (and the group is called {\sl parabolic}). In the hyperbolic group case, $h(\D)$ is symmetric with respect to the line of symmetry of the vertical strip, and ``quasi-symmetric'' with respect to any vertical line contained in the strip, and, in fact, the orbits of the group converge non-tangentially to the Denjoy-Wolff point. While, in the parabolic case, $h(\D)$ is highly non-symmetric with respect to any line contained in the half-plane and the orbits of the group converge tangentially to the Denjoy-Wolff point.

 For the general case, we will show that non-tangential convergence is equivalent to the image of the Koenigs function being ``quasi-symmetric'' about a vertical line. Suppose $\Omega \subsetneq \C$ is a domain starlike at infinity and $p \in \C$. Then for $t \geq 0$ define
\[
\delta^+_{\Omega,p}(t):=\min\{t,\ \inf\{|z-(p+it)|: z\in \partial \Omega, \Re z\geq \Re p\}\}, 
\]
and 
\[
\delta^-_{\Omega,p}(t):=\min\{t,\ \inf\{|z-(p+it)|: z\in \partial\Omega, \Re z\leq \Re p\}\}. 
\]
Then, the first main result we prove is  the following: 
\begin{theorem}\label{main3}
Let $(\phi_t)$ be a non-elliptic semigroup in $\D$ with Denjoy-Wolff point $\tau\in\partial \D$ and  Koenigs function $h$ and let $\Omega:=h(\D)$. Suppose that $\{t_n\}$ is a sequence converging to $+\infty$. Then 
\begin{enumerate}
\item the sequence $\{\phi_{t_n}(z)\}$ converges non-tangentially to $\tau$ as $n\to \infty$ for some---and hence any---$z\in \D$ if and only if for some---and hence any---$p\in \Omega$  there exist $0<c<C$ such that for all $n\in \N$
\[
c\delta^+_{\Omega,p}(t_n)\leq \delta^-_{\Omega,p}(t_n)\leq C\delta^+_{\Omega,p}(t_n).
\]
\item  $\lim_{n\to \infty}{\sf Arg}(1-\overline{\tau}\phi_{t_n}(z))=\frac{\pi}{2}$ (in particular,  $\{\phi_{t_n}(z)\}$ converges tangentially to $\tau$ as $n\to \infty$) for some---and hence any---$z\in \D$ if and only if for some---and hence any---$p\in \Omega$, 
\[
\lim_{n\to +\infty}\frac{\delta^+_{\Omega,p}(t_n)}{\delta^-_{\Omega,p}(t_n)}=0,
\]
while, $\lim_{n\to \infty}{\sf Arg}(1-\overline{\tau}\phi_{t_n}(z))=-\frac{\pi}{2}$ (in particular, $\{\phi_{t_n}(z)\}$ converges tangentially to $\tau$ as $n\to \infty$) for some---and hence any---$z\in \D$ if and only if for some---and hence any---$p\in \Omega$, 
\[
\lim_{n\to +\infty}\frac{\delta^+_{\Omega,p}(t_n)}{\delta^-_{\Omega,p}(t_n)}=+\infty.
\]
\end{enumerate}
\end{theorem}

The proof of this result is very involved, and it is based almost entirely on Gromov's theory of negatively curved metric spaces. In particular, let $k_\Omega$ denote the hyperbolic distance on $\Omega$. When $0\not\in\Omega$ and $it\in \Omega$ for all $t>0$, we show that the $2$-Lipschitz curve 
\begin{align}\label{eq:sigma}
\sigma:[1,+\infty)\ni t\mapsto  \frac{\delta^+_{\Omega,0}(t)-\delta^-_{\Omega,0}(t)}{2}+it
\end{align}
can be parametrized to be a quasi-geodesic in $(\Omega,k_\Omega)$  (see Section \ref{Sec:prel} for details on quasi-geodesics). Thus, by Gromov's shadowing lemma, $\sigma$ always stays within a finite hyperbolic distance from a geodesic ``converging to $\infty$.'' Theorem \ref{main3} then follows by noticing that non-tangential convergence is equivalent to  staying at finite hyperbolic distance from $\sigma$ (see Section  \ref{Sec:proof1} for details). 

Probably more  surprising than the previous result, we also prove that an orbit of a semigroup converges non-tangentially if and only if it can be parametrized to be a quasi-geodesic in $(\Omega, k_\Omega)$. More precisely, we show that the orbit $(\phi_t(z))$ converges non-tangentially if and only if   for every $0\leq t_1 \leq t_2$, the hyperbolic length of  the orbit of $(\phi_t(z))$ between $t_1$ and $t_2$ is, up to a fixed constant which does not depend on $t_1, t_2$, the hyperbolic distance between $\phi_{t_1}(z)$ and  $\phi_{t_2}(z)$.

Notice that,  in view of \cite{Bet, CDG, BCDG}, there exist examples of  semigroups   such  that the slope is an  interval $[a,b]$ with $-\pi/2< a<b< \pi/2$. Nonetheless, according to our result, the orbit can be parametrized to be a quasi-geodesic.

Collecting the results, the other two main results of the paper are the following theorems.

\begin{theorem}\label{main1}
Let $(\phi_t)$ be a non-elliptic semigroup in $\D$ with Denjoy-Wolff point $\tau\in\partial \D$ and  Koenigs function $h$ and let $\Omega:=h(\D)$. Then the following are equivalent:
\begin{enumerate}
\item for some---and hence any---$z\in \D$, the orbit $[0,+\infty)\ni t\mapsto \phi_t(z)$ converges non-tangentially to $\tau$ as $t\to +\infty$,
\item for some---and hence any---$z\in \D$, the curve $[0,+\infty)\ni t\mapsto \phi_t(z)$ can be parametrized to be a quasi-geodesic,
\item for some---and hence any---$p\in \Omega$ there exist $0<c<C$ such that for all $t\geq 0$,
\[
c\delta_{\Omega,p}^+(t)\leq \delta_{\Omega,p}^-(t)\leq C\delta_{\Omega,p}^+(t).
\]
\end{enumerate}
\end{theorem}

\begin{theorem}\label{main2}
Let $(\phi_t)$ be a non-elliptic semigroup in $\D$ with Denjoy-Wolff point $\tau\in\partial \D$ and  Koenigs function $h$ and let $\Omega:=h(\D)$. Then the following are equivalent:
\begin{enumerate}
\item $\lim_{t\to+\infty}{\sf Arg}(1-\phi_t(z))=\pi/2$ ({\sl respectively} $=-\pi/2$) for some---and hence any---$z\in \D$,   and, in particular, $[0,+\infty)\ni t\mapsto \phi_t(z)$ converges tangentially to $\tau$ as $t\to +\infty$,
\item  $\lim_{t\to +\infty}\frac{\delta_{\Omega,p}^+(t)}{\delta_{\Omega,p}^-(t)}=0$ ({\sl respect. } $\lim_{t\to +\infty}\frac{\delta_{\Omega,p}^+(t)}{\delta_{\Omega,p}^-(t)}=+\infty$).
\end{enumerate}
\end{theorem}

As we will show, Theorem \ref{main1} and Theorem \ref{main2} are consequences of  Theorem \ref{main3} and of its proof. 

Recall that a non-elliptic semigroup $(\phi_t)$ is {\sl hyperbolic} if  $h(\D)$ is contained in a vertical strip, it is {\sl parabolic  of positive hyperbolic step} if  $h(\D)$ is contained in a vertical half-plane but not in a vertical strip and {\sl parabolic  of zero hyperbolic step} otherwise. We mention that, although our proofs do not rely on previous results about dynamics of semigroups, it was already known (see \cite{CD, CDP}) that if $(\phi_t)$ is a hyperbolic semigroup then the trajectory $t\mapsto \phi_t(z)$ always converges non-tangentially to its Denjoy-Wolff point as $t\to +\infty$ for every $z\in \D$ and the slope is a single point which depends harmonically on $z$, while, if it is parabolic of positive hyperbolic step then $\phi_t(z)$ always converges tangentially to its Denjoy-Wolff point as $t\to +\infty$ for every $z\in \D$ and the slope is independent of $z$ (and it is either $\{\pi/2\}$ or $\{-\pi/2\}$).
 
 Therefore,  Theorem \ref{main2} gives the new information that every  orbit of a hyperbolic semigroup is a  quasi-geodesic, while, in the case of parabolic semigroups of positive hyperbolic step, the orbits are never  quasi-geodesics. 
 
 In the case of parabolic semigroups of zero hyperbolic step, all cases can happen.  In Section \ref{examples} we give some examples illustrating the possible behaviors. 
 
 The paper is organized as follows. In Section \ref{examples} we provide some examples of possible behavior of orbits. In Section \ref{Sec:prel} we state some preliminaries we need in this paper. In Section \ref{Quasi} we show that the curve $\sigma$ defined in Equation~\eqref{eq:sigma} can indeed be parametrized to be a quasi-geodesic and also estimate its hyperbolic distance to the vertical axis at $p$. Finally, in Section \ref{Sec:proof1} we prove the theorems.
 
 \section{Examples}\label{examples}
In this section we construct some examples  of parabolic semigroups of zero hyperbolic step illustrating possible cases. We define domains $\Omega$ starlike at infinity, and, if $h:\D \to \Omega$ is a Riemann map, the semigroup is given by $\phi_t(z):=h^{-1}(h(z)+it)$.

\begin{example}
The model domain $\Omega_1$ is defined by $\Omega_1:=\{\zeta \in \mathbb C : \Im(\zeta) > (\Re(\zeta))^2\}$ (see Figure \ref{Fig1}). 
 \begin{figure}[h!]


\begin{center}
\begin{tikzpicture}[scale=\ech]
\draw [domain=-1.8:1.8] plot (\x,{(\x)^2});
\draw [dotted] [domain= -2:-1.8] plot (\x,{(\x)^2});
\draw [dotted] [domain= 1.8:2] plot (\x,{(\x)^2});

\draw [dashed] (0,0.3) -- (0,4) ;

\draw (0,2) node[scale=0.65] {$\bullet$};
\draw (0.3,1.55) node[scale=0.75] {$it$};

\draw (-1.3,1.69) node[scale=0.65] {$\bullet$};
\draw (1.3,1.69) node[scale=0.65] {$\bullet$};

\draw[<->] (-1.25,1.68) -- (-0.05,1.98);

\draw (-0.7, 2.25) node[scale=0.65] {$\delta_{\Omega_1,0}^-(t)$};

\draw[<->] (1.25,1.68) -- (0.05,1.98);

\draw (0.8, 2.25) node[scale=0.65] {$\delta_{\Omega_1,0}^+(t)$};

\draw (-1.4,0.5) node{$\Omega_1$} ;

\end{tikzpicture}
\end{center}
  \caption{}\label{Fig1}
\end{figure}

Then $\Omega_1$ is symmetric with respect to the imaginary axis, $\delta^+_{\Omega,0}(t)=\delta^-_{\Omega,0}(t)$ for $t>0$ and $\gamma:[1,+\infty)\ni t  \mapsto it$ can be reparametrized as a geodesic  in $\Omega_1$. Hence, for every $z \in \mathbb D$, the semigroup $\phi_t(z)$ converges orthogonally to the Denjoy-Wolff point $\tau \in \partial \mathbb D$.
\end{example}

\begin{example}
The model domain $\Omega_2$ (see Figure \ref{Fig2}) is defined by
$$
\Omega_2:=\{\zeta \in \mathbb C:\Re(\zeta) > 0\} \cup \{\zeta \in \mathbb C :\Im(\zeta) > (\Re(\zeta))^2\}.
$$
\begin{figure}[h!]



\begin{center}
\begin{tikzpicture}[scale=\ech]
\draw [domain=-1.8:0] plot (\x,{(\x)^2});
\draw [dotted] [domain= -2:-1.8] plot (\x,{(\x)^2});

\draw (0,0) -- (0,-2) ;
\draw [dotted] (0,-2) -- (0,-2.5) ;

\draw [dashed] (0,0.3) -- (0,4) ;

\draw (0,2) node[scale=0.65] {$\bullet$};
\draw (0.25,2) node[scale=0.75] {$it$};

\draw (-1.3,1.69) node[scale=0.65] {$\bullet$};

\draw[<->] (-1.25,1.68) -- (-0.05,1.98);

\draw (-0.7, 2.25) node[scale=0.65] {$\delta_{\Omega_2,0}^-(t)$};

\draw (-0.8,-0.3) node{$\Omega_2$} ;

\end{tikzpicture}
\end{center}
    \caption{}\label{Fig2}
\end{figure}

Then for every $t > 4$, $\delta_{\Omega_2,0}^+(t) = t$ and  $\displaystyle \delta_{\Omega_2,0}^-(t) = \sqrt{t-\frac{1}{4}}$. Hence
$$
\lim_{t \rightarrow \infty} \frac{\delta_{\Omega_2,0}^+(t)}{\delta_{\Omega_2,0}^-(t)} = +\infty.
$$
It follows from Theorem~\ref{main2} that for every $z \in \mathbb D$, the semigroup $\phi_t(z)$ converges tangentially to the Denjoy-Wolff point $\tau \in \partial \mathbb D$.
\end{example}

\begin{example}
The model domain $\Omega_3$ (see Figure \ref{Fig3}) is defined by
$$
\Omega_3:=\Omega_2\cup_{n \geq 1}S_n,
$$
\begin{figure}[h!]



\begin{center}
\begin{tikzpicture}[scale=\ech]
\draw [domain=-0.5:0] plot (\x,{(\x)^2});
\draw (-0.5,0.25) -- (-0.5,-2) ;
\draw (-0.7,0.49) -- (-0.7,-2) ;


\draw [domain=-1:-0.7] plot (\x,{(\x)^2});
\draw (-1,1) -- (-1,-2) ;
\draw (-2.3,5.29) -- (-2.3,-2) ;

\draw [domain=-2.6:-2.3] plot (\x,{(\x)^2});

\draw [dotted] [domain=-2.7:-2.6] plot (\x,{(\x)^2});


\draw (0,0) -- (0,-2) ;
\draw [dotted] (0,-2) -- (0,-2.5) ;

\draw [dashed] (0,0.3) -- (0,4) ;

\draw (0,1.5) node[scale=0.65] {$\bullet$};
\draw (0.5,1.5) node[scale=0.75] {$it_1$};

\draw (-0.8,0.64) node[scale=0.65] {$\bullet$};

\draw (0,2.4) node[scale=0.65] {$\bullet$};
\draw (0.5,2.4) node[scale=0.75] {$is_1$};

\draw (-1,1) node[scale=0.65] {$\bullet$};

\draw[<->] (-0.97,1.02) -- (-0.05,2.38);

\draw (-1.25, 1.9) node[scale=0.65] {$\delta_{\Omega_3,0}^-(s_1)$};

\draw[<->] (-0.77,0.66) -- (-0.05,1.48);

\draw  (0.35, 0.8) node[scale=0.65] {$\delta_{\Omega_3,0}^-(t_1)$};

\draw (-3.5,-0.8) node {$S_{1}$} ;
\draw[->] (-3,-0.9) -- (-0.6,-0.6) ;

\draw (-3.5,0) node {$S_2$} ;
\draw[->] (-2.8,0) -- (-1.7,0.2) ;

\draw (-3.3,-1.8) node {$\Omega_3$};

\end{tikzpicture}
\end{center}
    \caption{}\label{Fig3}
\end{figure}

\noindent where for every $n \geq 1$, $S_n$ is a vertical strip $S_n:=\{\zeta \in \mathbb C: a_n < \Re(\zeta) < b_n < 0\}$. The sequences $(a_n)$ and $(b_n)$ are constructed inductively as follows.

First consider $t_1 > 4$. Then $\delta_{\Omega_2,0}^-(t_1) = |it_1 - \zeta_1|=\sqrt{t_1-\frac{1}{4}}$ with $\zeta_1= -\sqrt{t_1-\frac{1}{2}}+i(t_1-\frac{1}{2})$.

Let $b_1=-1 -\sqrt{t_1-\frac{1}{2}}$ and $\eta_1:=b_1+ib_1^2$. Select the unique $s_1 > t_1$ such that $|is_1 - \eta_1| = \frac{1}{2} s_1$. We can now choose $a_1 < b_1$ such that for every $\zeta \in \{z \in \C :\Re(z) \leq a_1,\Im(z) = (\Re(z))^2\}$, we have $|is_1 - \zeta| >  \frac{1}{2} s_1$. In particular we will have
$$
\delta_{\Omega_3,0}^-(s_1)=|is_1 - \eta_1| = \frac{1}{2} s_1 = \frac{1}{2} \delta_{\Omega_3,0}^+(s_1).
$$

Now, we may choose $t_2 >s_1$ and $\zeta_2 \in \partial \Omega_2$, with $\Re(\zeta_2) < a_1$, such that $\delta_{\Omega_2,0}^-(t_2) = |it_2 - \zeta_2|=\sqrt{t_2-\frac{1}{4}}$. Then let $b_2=-1 -\sqrt{t_2-\frac{1}{2}}$ and $\eta_2:=b_2+ib_2^2$. From $b_2$ and $\zeta_2$, we construct $a_2 < b_2$ and $s_2$, exactly as we constructed $a_1$ and $b_1$ from $b_1$ and $\zeta_1$. In particular we will have
$$
\delta_{\Omega_3,0}^-(s_2)=|is_2 - \eta_2| = \frac{1}{2} s_2 = \frac{1}{2} \delta_{\Omega_3,0}^+(s_2).
$$
The construction of sequences $(a_n)$ and $(b_n)$ is completed by induction.

By construction, $\lim_{n \rightarrow \infty}a_n = \lim_{n \rightarrow \infty}b_n = -\infty$ and $\lim_{n \rightarrow \infty}t_n = \lim_{n \rightarrow \infty}s_n = +\infty$.

Note that, for every $t\geq 0$, $\delta_{\Omega,0}^+(t)\geq \delta_{\Omega,0}^-(t)$, which, according to Theorem \ref{main3}, means that there are no subsequences of any orbit of $(\phi_t)$ converging to $\tau$ with slope $\pi/2$. On the other hand, for every $n \geq 1$, we have:
$$
\displaystyle \frac{\delta_{\Omega_3,0}^+(t_n)}{ \delta_{\Omega_3,0}^-(t_n)} = \displaystyle \frac{t_n}{\sqrt{t_n-\frac{1}{4}}} \to +\infty \quad \hbox{ as $n \rightarrow \infty$},
$$
which, again by Theorem \ref{main3}, implies that $\phi_{t_n}(z)\to \tau$ with slope $-\pi/2$. Finally, since
$$
\displaystyle \frac{\delta_{\Omega_3,0}^-(s_n)}{\delta_{\Omega_3,0}^+(s_n)} = \displaystyle \frac{s_n/2}{s_n} = \displaystyle \frac{1}{2},
$$
 for every $z \in \mathbb D$ the sequence $\{\phi_{s_n}(z)\}$ converges non-tangentially to $\tau$. In particular, the slope of $(\phi_t)$ is $[-\pi/2,\alpha]$ for some $-\pi/2<\alpha<\pi/2$.
\end{example}

\section{Preliminaries on hyperbolic and Euclidean geometry}\label{Sec:prel}

\subsection{Hyperbolic geometry of simply connected domains} Let $\Omega\subsetneq \C$ be a simply connected domain. Recall that  the hyperbolic metric $\kappa_\Omega$ is defined for $z\in \Omega$ and $v\in \C$ by
\[
\kappa_\Omega(z;v):=\frac{|v|}{f'(0)},
\]
where $f:\D\to \Omega$ is the Riemann map such that $f(0)=z$ and $f'(0)>0$. The hyperbolic distance between $z,w\in \Omega$ is defined as 
\[
k_\Omega(z,w):=\inf \int_0^1 \kappa_\Omega(\gamma(\tau);\gamma'(\tau))d\tau,
\]
where the infimum is taken over all piecewise $C^1$-smooth curves $\gamma:[0,1]\to \Omega$ such that $\gamma(0)=z$ and $\gamma(1)=w$.

A curve $\gamma:[a,b]\to \Omega$ is  rectifiable if 
\[
\ell_\Omega(\gamma;[a,b]):=\sup_{\mathcal P}\sum_{j=0}^N k_\Omega(\gamma(t_j), \gamma(t_{j+1}))<+\infty,
\]
where the supremum is taken over all partitions $\mathcal P$ of $[a,b]$ of type $a=t_0<t_1<\ldots<t_{N+1}=b$,  $N\in \N$.

The number $\ell_\Omega(\gamma;[a,b])$ is the hyperbolic length of $\gamma$ and, by definition, 
\[
\ell_\Omega(\gamma;[a,b])\geq k_\Omega(\gamma(a), \gamma(b)).
\] 

Every rectifiable curve can be reparametrized by hyperbolic arc length. 
If $\gamma$ is a Lipschitz curve then 
\[
\ell_\Omega(\gamma;[s,t])=\int_s^t \kappa_\Omega(\gamma(\tau);\gamma'(\tau))d\tau.
\]

\subsection{Geodesics and non-tangential convergence}
Let $-\infty\leq a<b\leq +\infty$. A smooth curve $\eta:(a,b)\to \Omega$ is a (unit speed) {\sl geodesic} if 
\[
t-s=k_\Omega(\eta(s), \eta(t))
\]
 for all $a<s<t<b$. 
 
Given $R>0$  and a geodesic $\eta:[0,+\infty)\to \Omega$, the {\sl hyperbolic sector around $\eta$ of amplitude $R$} is given by
\[
S_\Omega(\eta, R):=\{z\in \Omega: k_\Omega\left(z,\eta([0,+\infty))\right)<R\}.
\]

We can use hyperbolic sectors to detect non-tangential convergence (see for instance~\cite[Proposition 4.5]{BCDG}): 

\begin{proposition} \label{non-tg-sector-hyp}
Let $\Omega\subsetneq \C$ be a simply connected domain and let $f:\D\to \Omega$ be a Riemann map. 
\begin{enumerate}
\item Suppose  $\gamma:[0,+\infty)\to \Omega$ be a continuous curve such that $\lim_{t\to+\infty}k_\Omega(\gamma(0), \gamma(t))=+\infty$, then $f^{-1}(\gamma(t))$ converges non-tangentially to a point $\sigma\in \partial \D$ if and only if there exist $R>0$ and a geodesic $\eta:[0,+\infty)\to \Omega$ such that $\gamma(t)\in S_\Omega(\eta,R)$ for all $t$ sufficiently large.
\item Suppose $\{w_n\}\subset \Omega$ be a sequence such that $\lim_{n\to \infty}k_\Omega(w_0,w_n)=\infty$, then $w_n$ converges non-tangentially to a point $\sigma\in \partial \D$ if and only if there exist $R>0$ and a geodesic $\eta:[0,+\infty)\to \Omega$ such that $w_n \in S_\Omega(\eta,R)$ for all $n$ sufficiently large.
\end{enumerate}
\end{proposition}

\subsection{Quasi-geodesics}
 Given a general simply connected domain $\Omega \subsetneq \C$ it is essentially impossible to determine the geodesics in the hyperbolic metric. However, it is sometimes possible to find so-called quasi-geodesics which, by Gromov's shadowing lemma (also called Morse lemma, or the geodesic stability lemma), turn out to approximate geodesics. 

\begin{definition}
Let $-\infty<a <b\leq +\infty$. Let $\Omega\subsetneq \C$ be a simply connected domain and $\gamma:[a,b)\to \Omega$. Let $A\geq 1$, $B\geq 0$. We say that $\gamma$ is a {\sl $(A,B)$-quasi-geodesic} if   for all $a\leq s\leq t<b$,
\begin{align*}
\frac{1}{A} (t-s)-B\leq \ell_\Omega(\gamma;[s,t]) \leq A(t-s)+B.
\end{align*}
\end{definition}

For short, we say that $\gamma$ is a {\sl  quasi-geodesic} if there exist $A\geq 1, B\geq 0$ such that $\gamma$ is  a  $(A,B)$-quasi-geodesic.

By Gromov's shadowing lemma (see, {\sl e.g.},  \cite[Th\'eor\`eme 3.1, pag. 41]{Coo}) there exists $M>0$ (which depends only on $A, B$) such that if $\gamma:[0,+\infty)\to \Omega$ is a  $(A,B)$-quasi-geodesic  then there exists  a geodesic $\eta:[0,+\infty)\to \Omega$  such that $\eta(0)=\gamma(0)$ and for every $t\in [0,+\infty)$ 
\begin{equation}\label{Eq:shadows}
k_\Omega(\gamma(t), \eta([0,+\infty)))<M, \quad k_\Omega(\eta(t), \gamma([0,+\infty)))<M.
\end{equation}

\begin{remark}\label{quasi-geo-impl-nontg}
Let $\Omega\subsetneq \C$ be a simply connected domain and let $f:\D\to \Omega$ be a Riemann map. By the previous argument and Proposition \ref{non-tg-sector-hyp} it follows that if $\gamma:[0,+\infty)\to \Omega$ is a  quasi-geodesic  then $f^{-1}(\gamma(t))$ converges non-tangentially to a point $\sigma\in \partial \D$ as $t\to+\infty$.
\end{remark}

From the previous discussion, we have the following result which allows to detect quasi-geodesics:

\begin{proposition}\label{prop:reparam} Suppose that $\Omega\subsetneq \C$ is a simply connected domain and $\gamma:[0,+\infty)\to \Omega$ is a Lipschitz curve. If there exists $A \geq 1$ and $B\geq 0$ such that 
\[
\ell_\Omega(\gamma;[s,t])  \leq Ak_\Omega(\gamma(s), \gamma(t))+B 
\]
for all $0 \leq s \leq t$, then $\gamma$ can be reparametrized to be a $(A,B)$-quasi-geodesic. 
\end{proposition}

\subsection{Estimates on the hyperbolic distance}

As customary, for $p\in \Omega$ we let
\[
\delta_\Omega(p)=\inf\{|z-p|: z\in \C\setminus \Omega\}.
\]
In this paper we will use the following estimates for the hyperbolic metric and distance (see \cite[Section 3]{BCDG} for details):

\begin{theorem}[Distance Lemma]\label{Thm:Distance-Lemma} Let $\Omega\subsetneq \C$ be a simply connected domain. Then for every $z\in \Omega$ and $v\in \C$,
\[
 \frac{|v|}{4\delta_\Omega(z)}\leq \kappa_\Omega(z;v)\leq \frac{|v|}{\delta_\Omega(z)}.
\]
Moreover, for every $w_1, w_2\in \Omega$,
\[
 \frac{1}{4} \log \left(1+\frac{|w_1-w_2|}{\min\{\delta_\Omega(w_1), \delta_\Omega(w_2)\}} \right)\leq k_\Omega(w_1,w_2)\leq \int_{\Gamma}\frac{|dw|}{\delta_\Omega(w)},
\]
where $\Gamma$ is any absolutely continuous curve in $\Omega$ joining $w_1$ to $w_2$.
\end{theorem}

Note that Theorem \ref{Thm:Distance-Lemma} implies immediately that for all $z,w\in \Omega$, 
\begin{equation}\label{lem:dist_lower_bd}
k_\Omega(z,w) \geq \sup_{\zeta \in \C \setminus \Omega} \frac{1}{4}\abs{ \log \frac{ \abs{z-\zeta}}{\abs{w-\zeta}} }.
\end{equation}

\subsection{Euclidean geometry of domains starlike at infinity}\

A simply connected domain $\Omega\subsetneq \C$ is {\sl starlike at infinity} if $\Omega+it\subseteq \Omega$ for all $t\geq 0$. 

Let $\Omega$ be a simply connected domain which is starlike at infinity and  $p\in \C$. For $t>0$, let
\[
\tilde\delta_{\Omega,p}^+(t):=\inf\{|z-(p+it)|: \Re z\geq \Re p, z\in \C\setminus \Omega\},
\] 
\[
\tilde\delta_{\Omega,p}^-(t):=\inf\{|z-(p+it)|: \Re z\leq \Re p, z\in \C\setminus \Omega\}.
\] 
Note that, if $p+it\in \C\setminus \Omega$ then $\tilde\delta_{\Omega,p}^+(t)=\tilde\delta_{\Omega,p}^-(t)=0$. While, for $p\in \Omega$ and $t>0$, $\delta_\Omega(p+it)=\min\{\tilde\delta_{\Omega,p}^+(t), \tilde\delta_{\Omega,p}^-(t)\}$. 

Moreover, for $t>0$ we let
\[
\delta^+_{\Omega, p}(t):=\min\{\tilde\delta_{\Omega,p}^+(t), t\},\quad \delta^-_{\Omega,p}(t):=\min\{\tilde\delta_{\Omega,p}^-(t), t\}.
\]
Note that, since $\Omega$ is starlike at infinity, then $(0,+\infty)\ni t\mapsto \delta_{\Omega,p}^{\pm}(t)$ is non-decreasing. 

Simple geometric considerations allow to prove the following lemma:

\begin{lemma}\label{Lem:same-delta}
Let $\Omega$ be a simply connected domain starlike at infinity. For all $p, q\in \Omega$ there exist $0<c<C$  such that for all $t>0$ 
\[
c \delta^{\pm}_{\Omega,p}(t)\leq \delta^{\pm}_{\Omega,q}(t)\leq C\delta^{\pm}_{\Omega,p}(t).
\]
\end{lemma}

\section{Quasi-geodesics in starlike at infinity domains}\label{Quasi}

The aim of this section is to construct a quasi-geodesic in a domain $\Omega\subsetneq \C$ starlike at infinity which converges in the Carath\'eodory topology to ``$+\infty$'' and to get useful estimates on the hyperbolic distance from this curve to a vertical axis.

In all this section, we assume that $\Omega\subset \C$ is a domain starlike at infinity such that $0\not\in \Omega$ and $it\in \Omega$ for all $t>0$. 

We define $\sigma:[1,+\infty)\to \Omega$ by 
\begin{equation}\label{Eq:defino-sigma}
\sigma(t) := \frac{ \delta_{\Omega,0}^+(t)-\delta_{\Omega,0}^-(t)}{2} + it.
\end{equation}

\begin{lemma}\label{Lem:sigma-2-Lipschitz}
The curve $\sigma$ is $2$-Lipschitz. In particular, $|\sigma'(t)|\leq 2$ for almost every $t\geq 1$.
\end{lemma}
\begin{proof}
For all $s,t\geq 1$, using the triangle inequality we have $\delta_{\Omega,0}^\pm(t)\leq |t-s|+\delta_{\Omega,0}^\pm(s)$ and $\delta_{\Omega,0}^\pm(t)\geq -|t-s|+\delta_{\Omega,0}^\pm(s)$. Therefore,
\[
|\delta^{\pm}_{\Omega,0}(t)-\delta_{\Omega,0}^\pm(s)|\leq |t-s|. 
\]
From this it follows immediately that $\sigma$ is $2$-Lipschitz.
\end{proof}

\subsection{The curve $\sigma$ is up to reparametrization a quasi-geodesic.}

The aim of this subsection is to prove the following result:

\begin{theorem}\label{thm:QG} The curve $[1,+\infty)\ni t\mapsto \sigma(t)$ can be reparametrized to be a  quasi-geodesic in $\Omega$.
\end{theorem}

The proof  is rather long and technical and requires many lemmas.

Let 
\begin{align*}
\omega(t) := \delta_{\Omega,0}^+(t)+\delta_{\Omega,0}^-(t).
\end{align*}

\begin{lemma}\label{Lem:Andy1} For $t \geq 1$
\begin{align*}
\delta_\Omega(\sigma(t)) \geq\frac{1}{2\sqrt{2}} \omega(t). 
\end{align*}
\end{lemma}

\begin{proof} Fix $t \geq 1$. First consider the case  $\delta_{\Omega,0}^+(t) \geq \delta_{\Omega,0}^-(t)$, which implies that $\Re \sigma(t)\geq 0$.  

If $z \in \partial \Omega$ and $\Re (z) > 0$, then 
\begin{align*}
 \abs{z-\sigma(t)} \geq \abs{z-it} - \abs{it-\sigma(t)} \geq \delta_{\Omega,0}^+(t) - \frac{ \delta_{\Omega,0}^+(t)-\delta_{\Omega,0}^-(t)}{2} = \frac{\omega(t)}{2}. 
\end{align*}

Now, for $z \in \C$ define
\begin{align*}
\|z\|_1 = \abs{ \Re z}+\abs{\Im z}.
\end{align*}
Then 
\begin{align*}
\abs{z} \leq \|z\|_1  \leq \sqrt{2}\abs{z}
\end{align*}
 If $z \in \partial \Omega$ and $\Re z \leq 0$, then 
\begin{align*}
 \abs{z-\sigma(t)} \geq \frac{1}{\sqrt{2}} \| z-\sigma(t)\|_1.
 \end{align*}
 Further, since $\Re z \leq 0 \leq \Re \sigma(t)$ we have
\begin{align*}
\| z-\sigma(t)\|_1 & = \abs{\Re z-\Re \sigma(t)}+\abs{\Im z-\Im \sigma(t)} = \Re \sigma(t)-\Re z + \abs{\Im z - t} \\
& =  \Re \sigma(t) + \|z-it\|_1 \geq \Re \sigma(t) + \abs{z-it} \geq  \Re \sigma(t)+\delta_{\Omega,0}^-(t) \\
& = \frac{ \delta_{\Omega,0}^+(t)-\delta_{\Omega,0}^-(t)}{2} + \delta_{\Omega,0}^-(t) = \frac{1}{2} \omega(t).
\end{align*}
Hence 
\begin{align*}
 \abs{z-\sigma(t)} \geq \frac{1}{2\sqrt{2}} \omega(t).
 \end{align*}

The case when  $\delta_{\Omega,0}^+(t) \leq \delta_{\Omega,0}^-(t)$ is similar. 
\end{proof}

As a direct consequence of the previous lemma, Lemma \ref{Lem:sigma-2-Lipschitz} and Theorem \ref{Thm:Distance-Lemma}, we have:

\begin{lemma}\label{lem:length_upper_bd_1} If $1 \leq a < b < \infty$, then 
\begin{align*}
\ell_\Omega\left(\sigma; [a,b]\right) \leq 4\sqrt{2} \int_a^b \frac{1}{\omega(t)} dt.
\end{align*}
\end{lemma}

We can now prove Theorem \ref{thm:QG} in a simple case. 

\begin{proposition}\label{Prop:easy-case1} Suppose that there exist $\alpha,T_0> 0$ such that 
\begin{align*}
\omega(t) \geq \alpha t 
\end{align*}
for all $t \geq T_0$. Then $\sigma$ can be reparametrized to be a  quasi-geodesic in $\Omega$. 
\end{proposition}

\begin{proof} We have $\delta^\pm_{\Omega,0}(t)\leq t$ for all $t\geq 1$, hence, 
\[
\frac{|\delta^+_{\Omega,0}(t)-\delta^-_{\Omega,0}(t)|}{t}\leq \frac{|\delta^+_{\Omega,0}(t)|}{t}+\frac{|\delta^-_{\Omega,0}(t)|}{t}\leq 2.
\]
Therefore, for all $t\geq 1$,
\begin{align*}
t \leq \abs{\sigma(t)} \leq 2 t.
\end{align*}
So, by \eqref{lem:dist_lower_bd}, for all $1\leq a\leq b$, 
\begin{equation}\label{Eq:stima-sosotto}
 k_\Omega(\sigma(a), \sigma(b)) \geq \frac{1}{4}\abs{ \log \frac{\abs{\sigma(b)}}{\abs{\sigma(a)}} }  \geq \frac{1}{4}  \log\frac{b}{a}  -\frac{1}{4}\log 2.
 \end{equation}
 On the other hand, if $T_0 \leq a \leq b$, then by Lemma \ref{lem:length_upper_bd_1},
 \begin{equation}\label{Eq:formaggio}
 \ell_\Omega(\sigma; [a,b]) \leq 4\sqrt{2} \int_a^b \frac{dt}{\omega(t)} \leq \frac{4\sqrt{2}}{\alpha} \int_a^b \frac{dt}{t} = \frac{4\sqrt{2}}{\alpha}\log\frac{b}{a}.
 \end{equation}
From this last  inequality, \eqref{Eq:stima-sosotto}, and Proposition~\ref{prop:reparam} it follows at once that $\sigma$ can be reparametrized to be a quasi-geodesic in $\Omega$. 
\end{proof}

\begin{remark}\label{Rem:stima-buona-per-estim}
For future reference, we make the following observations. If there exist $\alpha,T_0> 0$ such that 
\begin{align*}
\omega(t) \geq \alpha t 
\end{align*}
for all $t \geq T_0$, then
\begin{enumerate}  
\item by the same token we obtained \eqref{Eq:stima-sosotto}, we have
\[
\max\{k_\Omega(ia, \sigma(b)), k_\Omega(\sigma(a), ib)\}\geq \frac{1}{4}  \log\frac{b}{a}  -\frac{1}{4}\log 2.
\]
Hence, by \eqref{Eq:formaggio}, there exist constants $A, B>0$ such that for every $T_0 \leq a\leq b$ we have
\begin{equation}\label{Eq:stima-formaggio1}
k_\Omega(\sigma(a),\sigma(b))\leq \ell_\Omega(\sigma;[a,b])\leq A\min\{k_\Omega(ia, \sigma(b)), k_\Omega(\sigma(a), ib)\}+B.
\end{equation}
\item Also, again arguing as in \eqref{Eq:stima-sosotto}, we have
\begin{equation}\label{Eq:Herve-does-not-like-formaggio}
 \int_a^b \frac{dt}{\omega(t)} \leq \frac{4}{\alpha}k_\Omega(ia, ib).
\end{equation}
\end{enumerate}

\end{remark}

Now we make the following assumption:\\

\noindent \textbf{Assumption:} there does {\sl not} exist $\alpha,T_0> 0$ such that $\omega(t) \geq \alpha t $
for all $t \geq T_0$. \\

Assuming this condition,  there exists $T_0>0$ such that $\omega(T_0)<T_0$. In particular,  $\max\{\delta_{\Omega,0}^+(T_0),\delta_{\Omega,0}^-(T_0)\} < T_0$. Hence, for every $t \geq T_0$ we have 
\[
\delta_{\Omega,0}^{\pm}(t) \leq t-T_0+ \delta_{\Omega,0}^{\pm}(T_0) < t-T_0+T_0=t.
\]
Therefore, for every $t\geq T_0$,
\begin{equation}\label{Eq:delta-min-t}
\delta_{\Omega,0}^{\pm}(t) < t.
\end{equation}

\subsubsection{Step 1: constructing sequences} Fix $a,b \in [T_0,\infty)$ with $a < b$.   We define a sequence of positive numbers $\{t_n\}$
\begin{align*}
a=t_0 < t_1<t_2< \dots 
\end{align*}
and complex numbers $\{z_n^\pm\}\subset \C\setminus \Omega$ such that 
\begin{align*}
y_n:=\Im z_n^+=\Im  z_n^-
\end{align*}
 with the following properties
\begin{enumerate}
 \item $\Re z_n^- < 0 < \Re z_n^+$, 
 \item $|\Re z_n^\pm|\leq \delta_{\Omega,0}^\pm(t_n)$,
 \item $y_n\leq t_n$,
 \item $\max\{ |\sigma(t_{n})-z_{n}^+|,|\sigma(t_{n})-z_{n}^-|\} \leq 2\omega(t_n)$,
 \item $\min\{|\sigma(t_n)-z_{n-1}^+|,|\sigma(t_n)-z_{n-1}^-|\}=6\omega(t_n)$ 
\end{enumerate}
for all $n \geq 1$. 

\medskip

 We first define $t_n$. We set $t_0:=a$. Assuming $t_{0}, \dots, t_{n-1}$ and $z_0^\pm, \dots, z_{n-1}^\pm$ have been selected satisfying properties (1) to  (5) above, we define $t_n$ as follows: 
 \begin{equation}\label{Eq:def-tn-ricors}
t_{n} := \max \left\{ t \geq t_{n-1} : \omega(s) \geq \frac{1}{6} \min\{ \abs{\sigma(s)-z_{n-1}^+},\abs{\sigma(s)-z_{n-1}^-} \} \text{ for all } s \in [t_{n-1}, t]\right\}.
\end{equation}

Note that if $n \geq 1$, then $t_{n-1} < t_{n} < +\infty$. Indeed, by Property (4) 
\begin{align*}
\min\left\{ \abs{\sigma(t_{n-1})-z_{n-1}^+},\abs{\sigma(t_{n-1})-z_{n-1}^-}\right\} \leq 2\omega(t_{n-1}),
\end{align*}
and  the continuity of $t \rightarrow \omega(t)$ implies that $t_{n} > t_{n-1}$.  

Further, since we are assuming that there does not exist $\alpha,T_0> 0$ such that $\omega(t) \geq \alpha t $
for all $t \geq T_0$, we must have $t_{n} < +\infty$.

Now we define $z_n^+$ and $z_n^-$. Assuming $t_{0}, \dots, t_n$ and $z_0^+,z_0^-, \dots, z_{n-1}^+, z_{n-1}^-$ have been selected we define $z_{n}^+, z_{n}^- \in \C \setminus \Omega$ as follows: first pick $a_{n}, b_{n} \in \C \setminus \Omega$ such that 
\begin{align*}
& { \rm Re}(a_{n})  \leq 0 \leq { \rm Re}(b_{n}), \\
& \abs{a_{n} - it_{n}}  = \delta_{\Omega,0}^-(t_{n}), \text{ and} \\
& \abs{b_{n} - it_{n}}  = \delta_{\Omega,0}^+(t_{n}).
\end{align*}
Since $t_n \geq T_0$, by \eqref{Eq:delta-min-t} we have 
$\Re(a_{n})  < 0 < \Re(b_{n})$. Then let 
\begin{align*}
y_{n} := \min \{ { \rm Im}(a_{n}), { \rm Im}(b_{n}) \}.
\end{align*}
Since $\Omega$ is starlike at infinity, $\max \{ { \rm Im}(a_{n}), { \rm Im}(b_{n}) \}\leq t_n$, hence $y_{n} \leq t_{n}$. Then define
\begin{align*}
z_{n}^+ := \Re(b_{n}) + i y_{n}, \quad z_{n}^- := \Re(a_{n}) + i y_{n}.
\end{align*}
Note that by construction, $\Re z_n^+\leq |b_n-it_n|= \delta^+_{\Omega,0}(t_n)$, and similarly  $|\Re z_n^-|\leq  \delta^-_{\Omega,0}(t_n)$.

Moreover,
\begin{equation}\label{Eq:estima-bff65}
\max\{|it_{n}-z_{n}^+|,|it_{n}-z_{n}^-|\} \leq \omega(t_{n}).
\end{equation}
Indeed, assume that $y_n=\Im b_n$ (a similar argument works in case $y_n=\Im a_n$). Hence, $|it_{n}-z_{n}^+|=\delta_{\Omega,0}^+(t_n)$, while
\begin{equation*}
\begin{split}
|it_{n}-z_{n}^-|&\leq |it_n-a_n|+|a_n-(\Re a_n+iy_n)|=\delta_{\Omega,0}^-(t_n)+(\Im a_n-y_n)\\&\leq \delta_{\Omega,0}^-(t_n)+(t_n-y_n)\leq \delta_{\Omega,0}^-(t_n)+|it_n-b_n|=\delta_{\Omega,0}^-(t_n)+\delta_{\Omega,0}^+(t_n)=\omega(t_n).
\end{split}
\end{equation*}

Also, clearly $|\sigma(t_{n})-it_{n}| \leq \omega(t_{n})$. This last inequality, together with \eqref{Eq:estima-bff65}, implies
\[
 |\sigma(t_{n})-z_{n}^\pm|\leq |\sigma(t_n)-it_n|+|it_n-z_n^\pm|\leq 2\omega(t_{n}).
\]
The construction is completed.

\subsubsection{Step 2: key estimates} We now establish key estimates on the sequences constructed in the previous step.

\begin{lemma}\label{lem:key_est_1} For $n \geq 1$ we have
\begin{align*}
3\omega(t_{n}) \leq y_{n}-t_{n-1} \leq \min\{ t_{n}-t_{n-1}, y_{n}-y_{n-1}\}.
\end{align*}
In particular, 
\begin{align*}
t_0 < y_1 \leq t_1 < y_2 \leq t_2 < \dots
\end{align*}
and $\lim_{n \rightarrow \infty} y_n = \infty$. 
\end{lemma}

\begin{proof} 
Fix $n\geq 1$. By property (5), 
\begin{align*}
\min\{|\sigma(t_n)-z_{n-1}^+|,|\sigma(t_n)-z_{n-1}^-|\}=6\omega(t_n)
\end{align*}
First assume that $|\sigma(t_n)-z_{n-1}^+|=6\omega(t_n)$. Then, by \eqref{Eq:estima-bff65} and taking into account that $\omega(t_n)\geq \omega(t_{n-1})$, we have 
\begin{align*}
 y_{n}-t_{n-1} & = |iy_{n}-it_{n-1}| \\
 & \geq  |\sigma(t_{n})-z_{n-1}^+| - |\sigma(t_{n})-iy_n| - | it_{n-1}- z_{n-1}^+| \\
 & \geq 6 \omega(t_{n}) -  2\omega(t_{n}) -  \omega(t_{n-1})  \geq \left(6 - 3 \right) \omega(t_{n})  = 3\omega(t_{n}).
\end{align*}
By property (3), $y_{n}-t_{n-1} \leq \min\{ t_{n}-t_{n-1}, y_{n}-y_{n-1}\}$. The case when $|\sigma(t_n)-z_{n-1}^-|=6\omega(t_n)$ is essentially the same.

Finally, the previous estimates show that  $\{y_n\}$ is an increasing sequence and
\[
0<3\omega(t_0)\leq 3\lim_{n\to\infty}\omega(t_n)\leq \lim_{n\to\infty}(y_n-y_{n-1}).
\]
Hence $\lim_{n \to \infty} y_n = \infty$. 
\end{proof}

As straightforward consequence of the previous lemma and taking into account that $\omega(t_n)\geq \omega(t_{n-1})$, we see that \begin{equation}\label{lem:key_est_3}
 \log \left( \frac{y_{n}-y_{n-1}}{\omega(t_{n-1})} \right) \geq \log 3>1
\end{equation}
for every $n \geq 1$. 

\begin{lemma}\label{lem:key_est_2} If $n \geq 1$ and $t \in [y_n, t_n]$, then 
\begin{align*}
\omega(t)\leq \omega(t_n) \leq 2 \omega(t).
\end{align*}
\end{lemma}

\begin{proof}
The first inequality follows from the fact that $\Omega$ is starlike at infinity.

Since $t_{n-1} < y_n \leq t_n$ it follows from \eqref{Eq:def-tn-ricors} and the fact that $\sigma$ is $2$-Lipschitz (see Lemma \ref{Lem:sigma-2-Lipschitz}) that  
   \begin{align*}
  \omega(t) 
  &\geq \frac{1}{6} \min\{ \abs{\sigma(t)-z_{n-1}^+}, \abs{\sigma(t)-z_{n-1}^-} \} \\
  & \geq \frac{1}{6} \min\{ \abs{\sigma(t_n)-z_{n-1}^+}, \abs{\sigma(t_n)-z_{n-1}^-} \} - \frac{1}{6} \abs{\sigma(t_n)-\sigma(t)}\\
&\geq \omega(t_n) -  2\frac{1}{6} (t_n-t) \geq \omega(t_n) -  \frac{1}{3}|it_n-iy_n|\geq \left(1-\frac{1}{3}\right)\omega(t_n) \geq \frac{1}{2} \omega(t_n),
 \end{align*}
and the proof is completed. 
\end{proof}

\subsubsection{Step 3: A lower bound on distance}

Define 
\begin{align*}
\delta_n := { \rm Re}(z_n^+) - { \rm Re}(z_n^-).
\end{align*}
By property (2) in the definition of the sequence $\{z_n^\pm\}$,
\begin{align*}
\delta_n \leq \delta^+_{\Omega,0}(t_n) + \delta^-_{\Omega,0}(t_n) = \omega(t_n).
\end{align*}

Fix $N \geq 0$ such that $y_{N} \leq b < y_{N+1}$ (recall that we fixed $a,b \in [T_0,\infty)$ with $a < b$). 

\begin{lemma}\label{lem:dist_lower_bd_sigma} Suppose $u \in \{ ia, \sigma(a)\}$ and $v \in \{ib, \sigma(b)\}$. If $N =0$, then 
\begin{align*}
 k_\Omega(u,v) \geq -\frac{1}{4} \log \left( 2\right)+ \frac{1}{4}\log \left(\max\left\{ 1,  \frac{b-y_0}{\omega(a)} \right\}\right).
\end{align*}
If $N \geq 1$, then 
\begin{align*}
 k_\Omega(u,v) \geq \frac{1}{4} \left(-\log 2 +  \log \left( \frac{ y_1-y_0}{\omega(a)} \right)+ \sum_{k=1}^{N-1} \log  \left(\frac{y_{k+1}-y_k}{\delta_k}\right)+\log \left( \max\left\{1, \frac{b-y_N}{\delta_N}\right\} \right) \right).
\end{align*}
\end{lemma}

\begin{proof} First suppose that $N=0$. If $b-y_0\leq \omega(a)$ 
 there is nothing to prove. So suppose that 
\begin{align*}
 \frac{b-y_0}{\omega(a)} \geq 1.
 \end{align*}
 
 By \eqref{Eq:estima-bff65} and property (4) in Step 1, 
 \begin{equation}\label{Eq:andy-11}
 \max \left\{ |u - z_0^+|, |u-z_0^-| \right\} \leq 2\omega(a).
 \end{equation}
 
 Next, since $|v-z_0^\pm|\geq |\Im v-\Im z_0^\pm|=b-y_0$, we have
 \begin{equation}\label{Eq:andy-12}
\min\left\{ |v - z_0^+|, |v-z_0^-| \right\} \geq b - y_0.
 \end{equation}
 
Putting together \eqref{lem:dist_lower_bd} with \eqref{Eq:andy-11} and \eqref{Eq:andy-12}, we have
\begin{equation*}
k_\Omega(u,v) \geq \frac{1}{4} \log \abs{\frac{v-z_0^+}{u-z_0^+}}\geq-\frac{1}{4} \log \left( 2\right)+ \frac{1}{4}\log \left( \frac{b-y_0}{\omega(a)} \right).
\end{equation*}

Next suppose that $N > 0$.  Let $\gamma:[0,T] \rightarrow \Omega$ be a unit speed geodesic  with $\gamma(0)=u$ and $\gamma(T)=v$. For $k=1,\dots, N$ define 
 \begin{align*}
  \tau_k := \min \{ t \geq 0 : { \rm Im}(\gamma(t)) = y_k \}. 
 \end{align*}
 Note that  $a < \tau_{1} < \tau_{2} < \dots < \tau_{N} < b$.

 Then, since $\Omega$ is starlike at infinity, 
 \begin{equation}\label{Eq:stay-in-bet-Andy}
 \Re( z_k^-)  < \Re(\gamma(\tau_k)) < \Re(z_k^+).
 \end{equation}
Also, since $|\gamma(\tau_{k+1}) - z_k^\pm|\geq |\Im \gamma(\tau_{k+1})-\Im z_k^\pm|=y_{k+1}-y_k$, we have 
\begin{equation}\label{Eq:andy121}
\min\left\{ |\gamma(\tau_{k+1}) - z_k^+|, |\gamma(\tau_{k+1})-z_k^-|\right\} \geq y_{k+1}-y_k.
\end{equation}
Moreover, by  \eqref{Eq:stay-in-bet-Andy} we have $|\gamma(\tau_{k}) - z_k^\pm|=|\Re \gamma(\tau_{k})-\Re z_k^\pm|\leq \delta_k$, hence
\begin{equation}\label{Eq:andy169}
\max\left\{ |\gamma(\tau_{k}) - z_k^+|, |\gamma(\tau_{k})-z_k^-| \right\} \leq \delta_k.
\end{equation}
 
Now, by \eqref{lem:dist_lower_bd}, \eqref{Eq:andy121} and \eqref{Eq:andy-11} we have
\begin{equation}\label{Eq:bound1Ngreat1}
k_\Omega(u, \gamma(\tau_{1})) \geq \frac{1}{4} \log \abs{\frac{\gamma(\tau_1)-z_0^+}{u-z_0^+}}\geq -\frac{1}{4} \log \left(2 \right)+\frac{1}{4} \log \left( \frac{ y_1-y_0}{\omega(a)} \right).
\end{equation}

For $k \geq 1$, \eqref{lem:dist_lower_bd}, \eqref{Eq:andy121} and \eqref{Eq:andy169} imply that 
\begin{equation}\label{Eq:bound1Ngreat2}
k_\Omega(\gamma(\tau_{k+1}), \gamma(\tau_{k})) \geq \frac{1}{4} \log \abs{\frac{\gamma(\tau_{k+1})-z_k^+}{\gamma(\tau_k)-z_k^+}}\geq \frac{1}{4} \log \left( \frac{ y_{k+1}-y_k}{\delta_k} \right).
\end{equation}

Finally, \eqref{lem:dist_lower_bd}, \eqref{Eq:andy169} implies that 
\begin{equation*}
k_\Omega(\gamma(\tau_{N}),v) \geq \frac{1}{4} \log \abs{\frac{v-z_N^+}{\gamma(\tau_N)-z_N^+}}\geq \frac{1}{4} \log \left(  \frac{b-y_N}{\delta_N} \right),
\end{equation*}
and hence 
\begin{equation}\label{Eq:bound1Ngreat3}
k_\Omega(\gamma(\tau_{N}), v)\geq \frac{1}{4}\log \left( \max\left\{1, \frac{b-y_N}{\delta_N}\right\} \right).
\end{equation}

Since $\gamma$ is a geodesic, we have 
\begin{align*}
 k_\Omega(u, v) =  k_\Omega(u, \gamma(\tau_{1})) + \sum_{k=1}^{N-1}  k_\Omega(\gamma(\tau_k), \gamma(\tau_{k+1})) +  k_\Omega(\gamma(\tau_N), v),
\end{align*}
The statement then follows from \eqref{Eq:bound1Ngreat1}, \eqref{Eq:bound1Ngreat2}, \eqref{Eq:bound1Ngreat3}.
\end{proof}

\subsubsection{Proof of Theorem~\ref{thm:QG}}\label{subsec:111}

By Lemma~\ref{lem:length_upper_bd_1} we have 
\begin{align*}
\ell_\Omega\left(\sigma; [a,b]\right) \leq 4\sqrt{2} \int_a^b \frac{1}{\omega(t)} dt.
\end{align*}
Hence to prove Theorem~\ref{thm:QG} it is enough to show that  $\int_a^b \omega(t)^{-1}dt$ is comparable to the lower bounds in Lemma~\ref{lem:dist_lower_bd_sigma}.

\begin{lemma}\label{Lem:Andy-11} If $T \in [a, y_1]$, then 
\begin{align*}
 \int_{a}^{T} \frac{dt}{\omega(t)} \leq 1+  6\log \left( \max\left\{ 1, \frac{T-y_0}{\omega(a)} \right\} \right)
 \end{align*}
 \end{lemma}
 
 \begin{proof}
Notice that by \eqref{Eq:estima-bff65}, $(a-y_0)=(t_0-y_0)\leq |it_0-z_0^\pm|\leq \omega(t_0)$, hence, 
\begin{align*}
y_0 \leq a \leq y_0 + \omega(a)
\end{align*}
and if $a \leq t$, then $\omega(a) \leq \omega(t)$. So
 \begin{align*}
 \int_a^{y_0+\omega(a)} \frac{dt}{\omega(t)} \leq  \int_a^{y_0+\omega(a)} \frac{dt}{\omega(a)} \leq 1.
 \end{align*}
 
 Now if $t \in [a, y_1]$, then by \eqref{Eq:def-tn-ricors},
 \begin{align*}
\omega(t) \geq \frac{1}{6} \min\{ \abs{\sigma(t)-z_0^+}, \abs{\sigma(t)-z_0^-} \} \geq \frac{1}{6} ( t-y_0).
\end{align*}
 So if $ T \geq y_0+\omega(a)$, then 
  \begin{align*}
 \int_{y_0+\omega(a)}^{T} \frac{dt}{\omega(t)} \leq 6\int_{y_0+\omega(a)}^{T} \frac{dt}{t-y_0} =  6 \log \left( \frac{T-y_0}{\omega(a)} \right). 
 \end{align*}
 \end{proof}

\begin{lemma}\label{Lem:Andy-12} For $k \geq 1$, 
\begin{align*}
 \int_{y_k}^{y_{k+1}} \frac{dt}{\omega(t)} \leq 8\log \left( \frac{y_{k+1}-y_k}{\omega(t_k)} \right).
 \end{align*}
 \end{lemma}
 
 \begin{proof}
 By Lemma~\ref{lem:key_est_1}, 
 \begin{align*}
 y_k + \omega(t_k) \leq  y_k + 3\omega(t_k) \leq y_{k+1}
 \end{align*}
 Further, by Lemma~\ref{lem:key_est_2}, if $t \in [y_k, t_k]$, then 
 \begin{align*}
 \omega(t) \geq \omega(t_k)/2
 \end{align*}
 and, since $\Omega$ is starlike at infinity, if $t \geq t_k$, then  $\omega(t) \geq \omega(t_k)$. Therefore,  $\omega(t) \geq \omega(t_k)/2$ when $t \geq y_k$. Thus
\begin{equation}\label{Eq:int-stima-1-Andy}
 \int_{y_k}^{y_{k}+\omega(t_k)} \frac{dt}{\omega(t)} dt \leq  \int_{y_k}^{y_{k}+\omega(t_k)} \frac{2dt}{\omega(t_k)} dt= 2.
 \end{equation}

Next consider $t \in [y_k+\omega(t_k),y_{k+1}]$. By \eqref{Eq:estima-bff65}, we have 
\[
t_k-y_k=|it_k-iy_k|\leq |it_k-z_k^\pm|\leq \omega(t_k).
\]
 Then $y_k+\omega(t_k) \geq t_k$. So $t\in [t_{k}, y_{k+1}]$ and $y_{k+1}\leq t_{k+1}$. Hence,  by \eqref{Eq:def-tn-ricors},
\begin{align*}
\omega(t) \geq  \frac{1}{6} \min\{ |\sigma(t)-z_k^+|, |\sigma(t)-z_k^-| \} \geq \frac{1}{6}  ( t-y_k).
\end{align*}
Therefore,
\begin{equation}\label{Eq:int-stima-2-Andy}
 \int_{y_k+\omega(t_k)}^{y_{k+1}} \frac{dt}{\omega(t)} dt \leq  6 \int_{y_k+\omega(t_k)}^{y_{k+1}} \frac{dt}{t-y_k}  = 6 \log \left( \frac{y_{k+1}-y_k}{\omega(t_k)} \right).
\end{equation}

Thus by \eqref{Eq:int-stima-1-Andy} and \eqref{Eq:int-stima-2-Andy} and \eqref{lem:key_est_3},
\begin{align*}
 \int_{y_k}^{y_{k+1}} \frac{dt}{\omega(t)} dt 
 & \leq 2 +6 \log \left( \frac{y_{k+1}-y_k}{\omega(t_k)} \right)  \leq 8 \log \left( \frac{y_{k+1}-y_k}{\omega(t_k)} \right), 
  \end{align*}
and we are done.
  \end{proof}

Repeating the proof of the previous lemma one can prove:

\begin{lemma}\label{Lem:Andy-13} If $N \geq 1$, then
\begin{align*}
\int_{y_N}^b \frac{dt}{\omega(t)} \leq 2 + 6\log \left( \max\left\{ 1, \frac{b-y_N}{\omega(t_N)} \right\} \right).
\end{align*}
\end{lemma}

Combining the estimates in the previous three lemmas we can estimate $\int_a^b \omega(t)^{-1}dt$. 

Recall that $a,b \in [T_0,\infty)$ with $a < b$ and  $N \geq 0$ is a natural number such that $y_{N} \leq b < y_{N+1}$.

If $N =0$, then Lemma \ref{Lem:Andy-11} implies
\begin{equation}\label{Eq:pandemonio1}
\ell_\Omega(\sigma;[a,b]) \leq 4\sqrt{2}+  24\sqrt{2}\log \left( \max\left\{ 1, \frac{b-y_0}{\omega(a)} \right\} \right),
\end{equation}
while if $N > 0$, then Lemma \ref{Lem:Andy-12} and Lemma \ref{Lem:Andy-13} imply
\begin{equation}\label{Eq:pandemonio2}
\begin{split}
\ell_\Omega(\sigma; [a,b])& \leq 12\sqrt{2}+ 24\sqrt{2}\log \left(\frac{y_1-y_0}{\omega(a)}  \right) \\&+ 32\sqrt{2}\sum_{k=1}^{N-1} \log \left( \frac{y_{k+1}-y_k}{\omega(t_k)} \right) + 24\sqrt{2}\log \left( \max\left\{ 1, \frac{b-y_N}{\omega(t_N)} \right\} \right).
\end{split}
\end{equation}

Then Lemma~\ref{lem:dist_lower_bd_sigma} and the fact that $\delta_k \leq \omega(t_k)$ imply that there exist $A>1$ and $B>0$ such that for every $T_0\leq a\leq b$,
\begin{align*}
\ell_\Omega(\sigma;[a,b]) \leq A k_\Omega(\sigma(a), \sigma(b)) + B.
\end{align*}
Now, since $\sigma([1,T_0])$ is compact, possibly taking a larger $B$, the previous estimate holds for every $1\leq a\leq b$, and Theorem \ref{thm:QG} is finally proved.

\begin{remark}
We also notice that by \eqref{Eq:pandemonio1}, \eqref{Eq:pandemonio2} and Lemma \ref{lem:dist_lower_bd_sigma}, there exist constants $A, B>0$ such that  for every $1\leq a\leq b$,
\begin{equation}\label{Eq:stima-formaggio2}
k_\Omega(\sigma(a),\sigma(b))\leq \ell_\Omega(\sigma;[a,b])\leq A\min\{k_\Omega(\sigma(a), ib), k_\Omega(\sigma(b), ia)\}+B.
\end{equation}
\end{remark}

As a consequence of the previous results, we have the following:

\begin{proposition}\label{Prop:beta-quasi-geodetica0}
Assume there exist $c,C>0$ such that $c\delta^-_{\Omega,0}(t)\leq \delta^+_{\Omega,0}(t)\leq C \delta^-_{\Omega,0}(t)$ for all $t\geq 1$. Then $\beta_i:[0,+\infty)\ni t\mapsto i+it$ can be reparametrized to be a quasi-geodesic.
\end{proposition}
\begin{proof}
Since $\delta_\Omega(it)=\min\{\delta_{\Omega,0}^-(t), \delta_{\Omega,0}^+(t)\}$ and $\delta^-_{\Omega,0}(t)$ is comparable to $\delta^+_{\Omega,0}(t)$, there exists $C'>1$ such that for every $t\geq 1$,
\[
\omega(t) \leq C'\delta_\Omega(it).
\]
In particular, by  Theorem \ref{Thm:Distance-Lemma}, we have for every $0\leq a\leq b$,
\[
\ell_\Omega(\beta_i;[a,b])\leq \int_a^b\frac{d\tau}{\delta_\Omega(i\tau)}\leq \frac{1}{C'}\int_a^b\frac{d\tau}{\omega(\tau)}.
\]
Therefore, in case there exist $\alpha, T_0>0$ such that $\omega(t)\geq \alpha t$ for all $t\geq T_0$, equation \eqref{Eq:Herve-does-not-like-formaggio} implies that $\beta_i$ can be reparametrized to be a quasi-geodesic.

On the other hand, if there exist no $\alpha, T_0>0$ such that $\omega(t)\geq \alpha t$ for all $t\geq T_0$, Lemmas \ref{lem:dist_lower_bd_sigma}, \ref{Lem:Andy-11}, \ref{Lem:Andy-12}, \ref{Lem:Andy-13} imply again that $\beta_i$ can be reparametrized to be a quasi-geodesic.
\end{proof}

\subsection{Estimates on the distance between $\sigma$ and the vertical axis}
For $t\geq 1$, let $s_t\in [1,+\infty)$ be such that 
\begin{equation}\label{Eq:defino-st}
k_\Omega(\sigma(s_t), it)=\min_{r\in [1,+\infty)}k_\Omega(\sigma(r), it).
\end{equation}

\begin{proposition}\label{Prop:estimate-Andy-qg}
There exist $\alpha>1, \beta > 0$ such that for every $t\geq 1$,
\[
k_\Omega(\sigma(t),it) \leq \alpha k_\Omega(\sigma(s_t), it)+\beta.
\]
\end{proposition}
\begin{proof}
Either by \eqref{Eq:stima-formaggio1}, or by \eqref{Eq:stima-formaggio2},  for $t\geq 1$, we have
\[
 k_\Omega(\sigma(t),\sigma(s_t))\leq Ak_\Omega(it, \sigma(s_t))+B.
\]
Therefore
\begin{align*}
k_\Omega(it, \sigma(t)) \leq k_\Omega(it, \sigma(s_t))+ k_\Omega(\sigma(s_t), \sigma(t)) \leq (A+1) k_\Omega(it, \sigma(s_t))+B,
\end{align*}
and we are done.
\end{proof}

\section{Proof of Theorem \ref{main3}, Theorem \ref{main1} and Theorem \ref{main2}} \label{Sec:proof1}

\subsection{Proof of Theorem \ref{main3}}
Let $(\phi_t), \tau, h, \Omega, \{t_n\}$ be as in Theorem \ref{main3}. 

We can suppose that $(\phi_t)$ is not a group of automorphisms of $\D$, for otherwise the result is clear. 

In this case, there exists $p\in \C$ such that $p\not\in \Omega$ and $p+it\in \Omega$ for all $t>0$. Up to a translation, we can assume $p=0$. In particular, this implies that $\tilde\delta_{\Omega,0}^\pm(t)=\delta_{\Omega,0}^\pm(t)$ for every $t>0$. 

\begin{lemma}\label{Lem:uno-vale-uno-sigh}
The sequence $\{\phi_{t_n}(h^{-1}(i))\}$ converges to $\tau$ as $n\to+\infty$ non-tangentially ({\sl respectively}, tangentially) if and only if for every $z\in \D$ the sequence $\{\phi_{t_n}(z)\}$ converges to $\tau$ as $n\to+\infty$ non-tangentially ({\sl respect.}, tangentially).
\end{lemma}
\begin{proof}
Since $k_\D(\phi_{t_n}(h^{-1}(i)), \phi_{t_n}(z)))\leq k_\D(h^{-1}(i), z)<+\infty$ for every $n\in \N$, it follows that $\phi_{t_n}(z)$ is contained in a fixed hyperbolic neighborhood of $\{\phi_{t_m}(h^{-1}(i)): m\in \N\}$ for all $n\in \N$. Therefore the result follows at once from the triangle inequality and from Proposition~\ref{non-tg-sector-hyp}.
\end{proof}

Let $\sigma$ be the curve defined in \eqref{Eq:defino-sigma}. 

\begin{lemma}\label{Lem:same-point-sig}
 $\lim_{t\to+\infty}h^{-1}(\sigma(t))=\tau$.
\end{lemma}
\begin{proof}
By Remark~\ref{quasi-geo-impl-nontg}, the limit $x:=\lim_{t\to+\infty}h^{-1}(\sigma(t))$ exists. Suppose for a contradiction that $x\neq \tau$. 

For $n\in \N$ consider the segments $\tilde C_n(s):=in+s\frac{ \delta_{\Omega,0}^+(n)-\delta_{\Omega,0}^-(n)}{2}$, $0\leq s\leq 1$. Note that $\tilde C_n\subset \Omega$ for all $n\in \N$. 

Let $C_n:=h^{-1}(\tilde C_n)$, $n\in \N$. Since $x\neq \tau$, the Euclidean diameter of $(C_n)$ is bounded from below by a constant $K>0$. 

Moreover,  for every sequence $\{z_n\}$ such that $z_n\in C_n$, it holds $\lim_{n\to +\infty}|h(z_n)|=\infty$.

Therefore, $(C_n)$ is a sequence of Koebe's arcs for $h$, contradicting the no Koebe arcs theorem (see, {\sl e.g.}, \cite[Corollary 9.1]{Pom}).
\end{proof}

\begin{corollary}\label{Cor:quando-converge}
 The sequence  $\{\phi_{t_n}(z_0)\}$ converges non-tangentially to $\tau$ as $n\to+\infty$ for all $z_0\in \D$ if and only if there exists $C>0$ such that for every $n\in\N$
\[
k_\Omega(it_n, \sigma([1,+\infty)))\leq C.
\]
Conversely, the sequence $\{\phi_{t_n}(z_0)\}$ converges tangentially to $\tau$ as $n\to+\infty$ for all $z_0\in \D$ if and only if for every $M>0$ there exists $n_M\geq 1$ such that for all $n\geq n_M$,
\[
k_\Omega(it_n, \sigma([1,+\infty)))>M.
\]
\end{corollary}
\begin{proof}
By Theorem \ref{thm:QG} the curve $\sigma$ can be parametrized to be a  quasi-geodesic in $\Omega$, hence by \eqref{Eq:shadows}, it is ``shadowed'' by a geodesic $\gamma$ in $\Omega$. The curve $h^{-1}(\gamma)$ is then a geodesic in $\D$ and by Lemma \ref{Lem:same-point-sig} it converges to $\tau$. Hence, by the triangle inequality and Proposition \ref{non-tg-sector-hyp}, the sequence $\{\phi_{t_n}(h^{-1}(i))\}$ converges non-tangentially to $\tau$ as $n\to +\infty$ if and only if it is contained in a hyperbolic neighborhood of $h^{-1}(\sigma[1,\infty))$. Since $h$ is an isometry for the hyperbolic distance, it follows that  $\{\phi_{t_n}(h^{-1}(i))\}$ converges non-tangentially to $\tau$ as $n\to +\infty$ if and only if there exists $C>0$ such that for every $n\in\N$
\[
k_\Omega(it_n, \sigma([1,+\infty)))\leq C.
\]
Conversely, since (again by Proposition \ref{non-tg-sector-hyp}) $\{\phi_{t_n}(h^{-1}(i))\}$ converges tangentially to $\tau$ as $n\to+\infty$ if and only if it is eventually outside any hyperbolic sector around $h^{-1}(\gamma)$, by the same token as before, we get that   $\{\phi_{t_n}(h^{-1}(i))\}$ converges tangentially to $\tau$ as $n\to+\infty$ if and only if for every $M>0$ there exists $n_M\in \N$ such that for all $n\geq n_M$,
\[
k_\Omega(it_n, \sigma([1,+\infty)))>M,
\]
and we are done.
 \end{proof}

Now, for $t\geq 1$ let $s_t$ be defined as in \eqref{Eq:defino-st}. Notice that
\[
k_\Omega(it, \sigma([1,+\infty)))=k_\Omega(it, \sigma(s_t))).
\]
Then,  by Proposition \ref{Prop:estimate-Andy-qg} and the Distance Lemma (see Theorem \ref{Thm:Distance-Lemma}), we have for all $t\geq 1$,
\[
k_\Omega(it, \sigma([1,+\infty)))\geq \frac{1}{\alpha} k_\Omega(it, \sigma(t))-\frac{\beta}{\alpha}\geq \frac{1}{4\alpha}\log \left( \frac{ |\delta_{\Omega,0}^+(t)- \delta_{\Omega,0}^-(t)|}{ 2\min \{\delta_\Omega(it), \delta_\Omega(\sigma(t))\} }\right)-\frac{\beta}{\alpha}.
\]
In other words, there exist $A, B>0$ such that for every $t\geq 1$,
\begin{equation}\label{Eq:stima-1FAK}
k_\Omega(it, \sigma([1,+\infty))\geq A \log \left( \frac{ |\delta_{\Omega,0}^+(t)- \delta_{\Omega,0}^-(t)|}{ 2\delta_\Omega(it) }\right)-B.
\end{equation}

Now, for $t\geq 1$ let $\eta_t:[0,1]\to \Omega$ be defined as
\[
\eta_t(r):=it+r \frac{ \delta_{\Omega,0}^+(t)-\delta_{\Omega,0}^-(t)}{2}.
\]
 For all $t\geq 1$ we have
\begin{equation}\label{Eq:stima-altoHerv}
k_\Omega(\sigma(s_t), it)\leq k_\Omega(it, \sigma(t))\leq \ell_\Omega(\eta_t;[0,1]).
\end{equation}
We compute $\ell_\Omega(\eta_t;[0,1])$. In order to do so, we need a lemma:
\begin{lemma}\label{Lem:estima-dist-eta-H}
For every $t\geq 1$ and for every $r\in [0,1]$ we have
\[
\delta_\Omega(\eta_t(r))\geq \delta_\Omega(it).
\]
\end{lemma}
\begin{proof}
Fix $t\geq 1$ and  assume that $\delta_{\Omega,0}^+(t)\geq \delta_{\Omega,0}^-(t)$ (the case $\delta_{\Omega,0}^+(t)\leq \delta_{\Omega,0}^-(t)$ is similar and we omit it). 

Fix $r\in [0,1]$. Notice that $\Re \eta_t(r)\geq 0$. Therefore, if $z\in \C\setminus \Omega$ and $\Re z\leq 0$, then 
\[
|\eta_t(r)-z|\geq |it-z|\geq \delta_{\Omega,0}^-(t)=\delta_\Omega(it).
\]
On the other hand, if $z\in \C\setminus \Omega$ and $\Re z> 0$, then $|it-z|\geq \delta_{\Omega,0}^+(t)$. Therefore,
\begin{equation*}\begin{split}
|\eta_t(r)-z|&\geq \inf_{|w-it|=\delta{\Omega,0}^+(t), \Re w>0}|\eta_t(r)-w|=\delta_{\Omega,0}^+(t)-\Re \eta_t(r)\\&
\geq \delta_{\Omega,0}^+(t)-\Re \sigma(t)=\frac{1}{2}\left(\delta^+_{\Omega,0}(t)+\delta^-_{\Omega,0}(t)\right)\geq \delta_\Omega(it),
\end{split}
\end{equation*}
and we are done.
\end{proof}

By Lemma \ref{Lem:estima-dist-eta-H} and the Distance Lemma (Theorem \ref{Thm:Distance-Lemma}), 
\begin{equation*}
\begin{split}
\ell_\Omega(\eta;[0,1])&=\int_0^1\kappa_\Omega(\eta(r);\eta'(r))dr\leq \frac{|\delta_{\Omega,0}^+(t)-\delta_{\Omega,0}^-(t)|}{2}\int_0^1\frac{dr}{\delta_\Omega(\eta(r))}\\ &\leq \frac{|\delta_{\Omega,0}^+(t)-\delta_{\Omega,0}^-(t)|}{2\delta_\Omega(it)}.
\end{split}
\end{equation*}
This latter equation together with \eqref{Eq:stima-altoHerv} and \eqref{Eq:stima-1FAK} implies that for every $t\geq 1$,
\begin{equation}\label{Eq:main-eq-dist}
 A\log \left( \frac{ |\delta_{\Omega,0}^+(t)- \delta_{\Omega,0}^-(t)|}{ 2\delta_\Omega(it) }\right)-B\leq k_\Omega(it, \sigma([1,+\infty))) \leq \frac{|\delta_{\Omega,0}^+(t)-\delta_{\Omega,0}^-(t)|}{2\delta_\Omega(it)}.
\end{equation}

The first part (the ``non-tangential part'') of Theorem \ref{main3} follows now directly from Corollary \ref{Cor:quando-converge} and \eqref{Eq:main-eq-dist}. 

Also, by the same token, we see that $\phi_{t_n}(z)\to \tau$ tangentially  if and only if $\frac{\delta^+_{\Omega,0}(t_n)}{\delta^-_{\Omega,0}(t_n)}$  converges either to $0$ or $+\infty$ as $n\to \infty$.

We are left to show that 
\begin{equation}\label{Eq:deltaoverdelta1}
\lim_{n\to \infty}\frac{\delta^+_{\Omega,0}(t_n)}{\delta^-_{\Omega,0}(t_n)}= +\infty
\end{equation}
 if and only if 
 \begin{equation}\label{Eq:deltaoverdelta2}
 \lim_{n\to\infty}{\sf Arg}(1-\overline\tau \phi_{t_n}(z))=-\frac{\pi}{2}.
 \end{equation}

To this aim, we extend $\sigma$ to all of $(0,\infty)$ in the obvious way:
\begin{align*}
\sigma(t) = \frac{\delta^+_{\Omega,0}(t)-\delta^-_{\Omega,0}(t)}{2}+it
\end{align*}
Since 
 $0\not\in \Omega$ and $it\in \Omega$ for all $t>0$,  $\lim_{t\to 0^+}\sigma(t)=0$. Then $\sigma((0,\infty))$ divides $\Omega$ into the connected domains 
 \begin{align*}
 U^+ = \{ x+iy \in \Omega : x > \Re \sigma(y)\}
 \end{align*}
 and 
  \begin{align*}
 U^- = \{ x+iy \in \Omega : x < \Re \sigma(y)\}.
 \end{align*}
 Hence,  $\Gamma:=\overline{h^{-1}(\sigma(0,+\infty))}$ divides $\D$ into two connected components $D^+:=h^{-1}(U^+)$ and $D^-:=h^{-1}(U^-)$.

Also, there exists $\tilde\tau\in \partial \D$, $\tilde\tau\neq \tau$ such that $\lim_{t\to 0^-}h^{-1}(\sigma(t))=\tilde\tau$ (see, {\sl e.g.}, \cite[Theorem 1, p. 37]{Go}). 

By Remark \ref{quasi-geo-impl-nontg} and Lemma \ref{Lem:same-point-sig}, $h^{-1}(\sigma(t))$ converges to $\tau$ non-tangentially as $t\to+\infty$. This implies  that $\Gamma$ is contained in the set 
\begin{align*}
 \{z\in \D: |{\sf Arg}(1-\overline{\tau}z)|\leq \theta\}\cup \{\tau,\tilde\tau\}
\end{align*}
for some  $\theta\in (0,\pi/2)$. Notice that the last set is an angular sector of amplitude $2\theta$ with vertex $\tau$ symmetric with respect to segment joining $-\tau$ with $\tau$. 

Since $h$ preserves orientation, it follows that $D^+$ contains all the sequences converging tangentially to $\tau$ with slope $\pi/2$ while $D^-$ contains all the sequences converging tangentially to $\tau$ with slope $-\pi/2$.

Therefore, if \eqref{Eq:deltaoverdelta1} holds, then  $it_n\in U^-$ for $n$ sufficiently big, hence, $\phi_{t_n}(z)\in D^-$ eventually and \eqref{Eq:deltaoverdelta2} holds. Conversely, if \eqref{Eq:deltaoverdelta2} holds then $\phi_{t_n}(z)\in D^-$ eventually, hence, $it_n\in U^-$ eventually and  \eqref{Eq:deltaoverdelta1} holds.

This concludes the proof of the theorem.

\subsection{Proof of Theorem \ref{main1}}

The part ``(1) if and only if (3)'' follows immediately from Theorem \ref{main3}. By Remark \ref{quasi-geo-impl-nontg} it is clear that (2) implies (1) in Theorem \ref{main1}. In order to end the proof, we show that (3) implies (2). 

We need to prove  that the orbit $[0,+\infty)\ni t\mapsto \phi_t(z)$ can be parametrized to be a  quasi-geodesic for every $z\in \D$. Since $h$ is an isometry between $k_\D$ and $k_\Omega$, the latter statement is equivalent to proving that, setting $p=h(z)$, the curve $\beta_p:[0,+\infty)\ni t\mapsto p+it$ can be parametrized to be a  quasi-geodesic in $\Omega$. 

As before, we can assume $0\not\in \Omega$ and $it\in \Omega$ for all $t>0$. Hence, by Proposition \ref{Prop:beta-quasi-geodetica0}, the curve $\beta_i$ can be parametrized to be a quasi-geodesic in $\Omega$. 

In order to complete the proof, we will prove the following:

\begin{lemma} For every $p \in \Omega$, there exists $A_p > 1$ and $B_p > 0$ such that 
\begin{equation*}
\ell_\Omega(\beta_p;[s,t])\leq A_pk_\Omega(\beta_p(s),\beta_p(t))+B_p,
\end{equation*}
for all $0\leq s\leq t$. Hence, by Proposition~\ref{prop:reparam}, $\beta_p$ can be parametrized to be a  quasi-geodesic in $\Omega$. 
\end{lemma}

\begin{proof} Fix $p \in \Omega$. By Proposition~\ref{Prop:beta-quasi-geodetica0},  there exists $A\geq 1$ and $B\geq 0$ such that 
\begin{equation}\label{Eq:quasi-g-FFF}
\ell_\Omega(\beta_i;[s,t])\leq Ak_\Omega(\beta_i(s),\beta_i(t))+B,
\end{equation}
for all $0\leq s\leq t$. Now, for $0\leq s\leq t$,
\begin{equation*}
\begin{split}
k_\Omega(i+is,i+it)&\leq k_\Omega(p+is,i+is)+k_\Omega(p+is,p+it)+k_\Omega(i+it,p+it)\\ &\leq k_\Omega(p+is,p+it)+2k_\Omega(p,i),
\end{split}
\end{equation*}
where the last inequality follows from the fact that $\Omega \ni z\mapsto z+it$ is a holomorphic self-map of $\Omega$. Therefore, there exists $B_1>0$ such that for all $s,t\geq 0$,
\begin{equation}\label{Eq:Fil01}
k_\Omega(i+is,i+it)\leq k_\Omega(p+is,p+it)+B_1.
\end{equation} 

By Lemma  \ref{Lem:same-delta} there exists $c>0$ such that $\delta_\Omega(i+it)\leq c \delta_\Omega(p+it)$ for all $t\geq 0$. 
Hence, by the Distance Lemma (Theorem \ref{Thm:Distance-Lemma}), for $0\leq s\leq t$,
\begin{equation*}
\begin{split}
\ell_\Omega(\beta_p;[s,t])&=\int_s^t\kappa_\Omega(\beta_p(r);\beta'_p(r))dr\leq \int_s^t\frac{dr}{\delta_\Omega(p+ir)}\\&\leq c \int_s^t\frac{dr}{\delta_\Omega(i+ir)}\leq 4c\int_s^t\kappa_\Omega(\beta_i(r);\beta'_i(r))dr=4c \ell_\Omega(\beta_i;[s,t]).
\end{split}
\end{equation*}
Therefore, by \eqref{Eq:quasi-g-FFF} and \eqref{Eq:Fil01}
\begin{equation*}
\begin{split}
\ell_\Omega(\beta_p;[s,t])&\leq 4c \ell_\Omega(\beta_i;[s,t])\leq 4cAk_\Omega(i+is,i+it)+4cB\\&\leq 4cAk_\Omega(p+is,p+it)+4cAB_1+4cB,
\end{split}
\end{equation*}
for all $0\leq s\leq t$.
\end{proof}

\subsection{Proof of Theorem \ref{main2}} It follows directly from Theorem \ref{main3}.


\begin{thebibliography}{99}

\bibitem{Abate} M. Abate, {\sl Iteration theory of holomorphic maps on taut manifolds}, Mediterranean Press, Rende, 1989.


\bibitem{BrAr} L. Arosio, F. Bracci, {\sl Canonical models for holomorphic iteration}. Trans. Amer. Math. Soc., 368, (2016), 3305--3339.

\bibitem{Berkson-Porta} E. Berkson, H. Porta, {\sl Semigroups of
holomorphic functions and composition operators}. Michigan
Math. J., 25 (1978), 101--115.

\bibitem{Bet}  D. Betsakos, {\sl On the asymptotic behavior of the trajectories of semigroups of holomorphic functions}. J. Geom. Anal., 26 (2016),  557--569.

\bibitem{BCD} F. Bracci, M. D. Contreras, S. D\'iaz-Madrigal, {\sl Evolution Families and the Loewner Equation I: the unit disc}. J. Reine Angew. Math.(Crelle's Journal), 672, (2012), 1-37. 

\bibitem{BCDG} F. Bracci, M. D. Contreras, S. D\'iaz-Madrigal, H. Gaussier, {\sl Non-tangential limits and the slope of trajectories of holomorphic semigroups of the unit disc}. ArXiv 1804.05553, 2018

\bibitem{CD} M. D. Contreras, S. D\'iaz-Madrigal, {\sl Analytic flows in the unit disk: angular derivatives and boundary
fixed points}. Pacific J. Math., 222 (2005), 253--286.

\bibitem{CDG} M. D. Contreras, S. D\'iaz-Madrigal, P. Gumenyuk, {\sl Slope problem for trajectories of holomorphic semigroups in the unit disk}. Comput. Methods Funct. Theory, 15 (2015), 117--124.

\bibitem{CDP} M. D. Contreras, S. D\'iaz-Madrigal,  Ch. Pommerenke, {\sl Some remarks on the Abel equation in the
unit disk}, J. London Math. Soc. (2), 75, (2007), 623--634.

\bibitem{Coo}  M. Coornaert, T. Delzant, A. Papadopoulos,  {\sl G\'eom\'etrie et  th\'eorie des groupes. Les groupes hyperboliques de Gromov. (Geometry and group theory. The hyperbolic groups of Gromov)}.  Lecture Notes in Mathematics Series 1441, Springer-Verlag (1990), 165 p. 

\bibitem{Cowen} C. C. Cowen, {\sl Iteration and the solution of functional equations for functions analytic in the unit disk}, Trans. Amer. Math. Soc., 265 (1981), 69--95.

\bibitem{EliShobook10} M.\,Elin and D.\,Shoikhet, {\sl Linearization models for complex dynamical systems}. Topics in univalent functions, functional equations and semigroup theory. Birkh\"auser Basel, 2010.

\bibitem{GdeL}  E. Ghys, P. de La Harpe, {\sl Sur les groupes hyperboliques d'apr\`es Mikhael Gromov}, Progress in Mathematics,
83, Birkh\"auser.

\bibitem{Go} G. M. Goluzin, G.M.: {\sl Geometric Theory of Functions of a Complex Variable}. American Mathematical Society, Providence, R.I. (1969). (translated from G. M. Goluzin, Geometrical theory of functions of a complex variable (Russian), Second edition, Izdat. ``Nauka'', Moscow, 1966)

\bibitem{KaMc1} S. Karlin and J.  McGregor, {\sl Embeddability of discrete time simple branching processes into continuous time branching processes}. Trans. Amer. Math. Soc. 132 (1968), 115-136.

\bibitem{KaMc2} S. Karlin and J.  McGregor, {\sl Embedding iterates of analytic functions with two fixed points into continuous groups}. Trans. Amer. Math. Soc. 132 (1968), 137-145.

\bibitem{Mit} G. Mitchison, {\sl Conformal growth of Arabidopsis leaves}. Journal of Theoretical Biology 408 (2016), 155-166.

\bibitem{Pom}  Ch. Pommerenke, {\sl Univalent functions}, Vandenhoeck and Ruprecht, G{\"o}ttingen, 1975

\bibitem{ReiShobook05} S. Reich and  D. Shoikhet, {\sl Nonlinear semigroups, fixed points, and geometry of domains in Banach spaces}, Imperial College Press, 2005.

\bibitem{Shb} D. Shoikhet, {\sl Semigroups in geometrical function theory}. Kluwer Academic Publishers, Dordrecht, 2001.


\bibitem{Wo} J. Wolff, {\sl L'\'equation diff\'erentielle $dz/dt=w(z)=$fonction holomprphe \`a partie r\'eelle positive dans un demi-plan}, Compos. Math. \textbf{6} (1939), 296--304.


\end{thebibliography}
\end{document}